\documentclass[preprint,10pt]{elsarticle}

\usepackage{mathrsfs,amsmath,amssymb}

\usepackage{mathrsfs,amsmath,amssymb}
\usepackage{mathrsfs,amsmath,amssymb}
\usepackage{amscd} 
\usepackage{relsize,amsmath} 

\makeatletter
\newcommand{\biggg}[1]{{\hbox{$\left#1\vbox to 20.5pt{}\right.\n@space$}}}
\newcommand{\Biggg}[1]{{\hbox{$\left#1\vbox to 23.5pt{}\right.\n@space$}}}
\newcommand{\bigggg}[1]{{\hbox{$\left#1\vbox to 26.5pt{}\right.\n@space$}}}
\newcommand{\Bigggg}[1]{{\hbox{$\left#1\vbox to 29.5pt{}\right.\n@space$}}}
\newcommand{\biggggg}[1]{{\hbox{$\left#1\vbox to 32.5pt{}\right.\n@space$}}}
\newcommand{\Biggggg}[1]{{\hbox{$\left#1\vbox to 35.5pt{}\right.\n@space$}}}
\newcommand{\bigggggg}[1]{{\hbox{$\left#1\vbox to 38.5pt{}\right.\n@space$}}}
\newcommand{\Bigggggg}[1]{{\hbox{$\left#1\vbox to 41.5pt{}\right.\n@space$}}}
\makeatother

\makeatletter
\@addtoreset{equation}{section}
\makeatother

\usepackage{braket}

\begin{document}

\newtheorem{thm}{Theorem}
\newtheorem{lem}[thm]{Lemma}
\newdefinition{rmk}{Remark}
\newproof{pf}{Proof}
\newproof{pot}{Proof of Theorem \ref{thm2}}

\begin{frontmatter}



\title{Asymptotic behavior of solutions toward the constant state 
to the Cauchy problem for the non-viscous diffusive dispersive conservation law}


\author[labe1
]{Natsumi Yoshida}
\ead{14v00067@gmail.com
}

\address[label]{Graduate Faculty of Interdisciplinary Research Faculty of Education, 
University of Yamanashi, Kofu, Yamanashi 400-8510, Japan.}
\address{}

\begin{abstract}
In this paper, we investigate the asymptotic behavior of solutions 
to the Cauchy problem 
for the scalar non-viscous diffusive dispersive conservation laws 
where the far field states are prescribed. 
We proved that 
the solution of the Cauchy problem tends toward the constant state 
as time goes to infinity.
\end{abstract}


\end{frontmatter}

%



\pagestyle{myheadings}
\thispagestyle{plain}
\markboth{N. YOSHIDA}{NON-VISCOUS DIFFUSIVE DISPERSIVE CONSERVATION LAWS}

\section{Introduction and main theorems}
In this paper, 
we consider the asymptotic behavior of solutions to the Cauchy problem 
for a one-dimensional scalar diffusive dispersive conservation laws 
without viscous flux 
\begin{eqnarray}
 \left\{\begin{array}{ll}
  \partial_tu +\partial_x \bigl( \, 
f(u) + \delta \, \partial_x^2u 
+ \nu \, \partial_x^3u 
 \, \bigr) =0
  \qquad &\big( \, t>0, \: x\in \mathbb{R} \, \big), \\[5pt]
 u(0, x) = u_0(x) \rightarrow \tilde{u} \qquad &( x \rightarrow \pm \infty ),
 \end{array}
 \right.\,
\end{eqnarray}
where, $u=u(t, x)$ is the unknown function of $t>0$ and $x\in \mathbb{R}$, 
the so-called conserved quantity, 
$\tilde{u}\in \mathbb{R}$ is the constant state, 
$$
f(u) + \delta \, \partial_x^2u 
+ \nu \, \partial_x^3u 
\quad \big( \, \delta \in \mathbb{R}, \; \nu \geq 0 \, \big)
$$ 
is the total flux 
(that is, the functions $f(u)$, 
$\delta \, \partial_x^2u$ 
and $\nu \, \partial_x^3u$ stand for 
the convective flux, 
dispersive one 
and diffusive one, respectively), 
$u_0$ is the initial data, 
and $u_{\pm } \in \mathbb{R}$ 
are the prescribed far field states. 
We suppose that $f$ is a smooth function. 
It is noted that, the equation in the problem (1.1) is the non-viscous case ($\mu=0$) for the 
following equation:
\begin{equation}
\partial_tu +\partial_x \bigl( \, 
f(u) - \mu \, \partial_xu + \delta \, \partial_x^2u 
+ \nu \, \partial_x^3u 
 \, \bigr) =0, 
\end{equation} 
where $\mu \, \partial_xu$ is viscous/diffusive flux 
($\mu$ is the so-called viscous coefficient or anti-diffusion coefficient). 
It should be noted that, in the case $\mu > 0$, $\delta =0$, $\nu = 0$, 
(1.2) becomes the viscous conservation law/generalized viscous Burgers equation, 
in the case $\mu = 0$, $\delta \in \mathbb{R}$, $\nu = 0$, 
the one does the Korteweg-de Vries equation as one of the dispersive conservation laws, 
in the case $\mu > 0$, $\delta \in \mathbb{R}$, $\nu = 0$, 
the one does the generalized Korteweg-de Vries-Burgers equation,
in the case $\mu > 0$, $\delta \in \mathbb{R}$, $\nu > 0 $, 
the one does the generalized Korteweg-de Vries-Burgers-Kuramoto equation or 
the derivative form of the Kuramoto-Sivashinsky equation. 
We also note that the Korteweg-de Vries equation
can be categorized as dispersive conservation laws,
and the Korteweg-de Vries-Burgers equation 
and the Korteweg-de Vries-Burgers-Kuramoto equation or 
the derivative form of the Kuramoto-Sivashinsky equation the 
diffusive dispersive conservation laws. 

There have been known the various of the stability results 
concerning with the conservation laws (see \cite{and-ego-lan-tes}, \cite{bona-rajopadhye-schonbek}, 
\cite{du-gu}, \cite{dua-zha}, \cite{ego-gru-tes}, \cite{ego-tes}, 
\cite{harabetian}, \cite{hashimoto-matsumura}, \cite{hattori-nishihara}, \cite{ilin-oleinik}, 
\cite{lax}, \cite{liu-matsumura-nishihara}, 
\cite{matsu-nishi1}, \cite{matsu-nishi2}, \cite{matsu-nishi3}, 
\cite{matsumura-yoshida}, \cite{matsumura-yoshida'}, 
\cite{nishi-raj}, 
\cite{rajopadhye}, 
\cite{wan-zhu}, 
\cite{yoshida1}, \cite{yoshida2}, \cite{yoshida3}, \cite{yoshida4}, \cite{yoshida5}, 
\cite{yoshida7}, \cite{yoshida8} 
and so on, cf. \cite{chh1}, \cite{chh2}, 
\cite{chh-ric}, \cite{gur-mac}, \cite{ilin-kalashnikov-oleinik}, \cite{jah-str-mul}, 
\cite{liep-rosh}, 
\cite{ma}, \cite{ma-pr-st}, \cite{smoller}, \cite{soc}). 
In partucular, for the Cauchy problem of (1.2), 
Duan-Fan-Kim-Xie \cite{dua-fan-kim-xie}, 
Ruan-Gao-Chen \cite{rua-gao-che}, 
showed some stabilities of the rarefaction waves and 
Yoshida \cite{yoshida6} showed the global stabilities of the constant state 
and the rarefaction wave. 

However, the any stabilities have not yet been known for more difficult non-viscous case, (1.1). 

Our main theorem is stated as follows.  

\medskip

\noindent
{\bf Theorem 1.1} (Main Theorem ){\bf .}\quad{\it
Assume 
the convective flux $f \in C^2(\mathbb{R})$ satisfy 
\begin{equation}
|\, f''(u) \,| \leq O(1)\, \big( \, 1 + |\, u \,|^q \, \big)
\quad \big( \, 0\leq q\leq5 \, \big) 
\end{equation}
and the initial data satisfy
$u_0-\tilde{u} \in L^2$ and
$\partial _xu_0 \in H^1$. 
Then the Cauchy problem {\rm(1.1)} has a 
unique global in time 
solution $u$ 
satisfying 
\begin{eqnarray*}
 \left\{\begin{array}{ll}
u-\tilde{u} \in C^0\bigl( \, [\, 0, \, \infty \, ) \, ; H^2 \, \bigr), \\[5pt]
\partial _x^2 u \in L^2\bigl( \, 0,\, \infty \, ; H^2 \, \bigr),
\end{array} 
\right.\,
\end{eqnarray*}
and the asymptotic behavior 
$$
\lim _{t \to \infty}
\bigg( \,
\sup_{x\in \mathbb{R}}
| \, u(t,x) - \tilde{u} \, | 
+\sup_{x\in \mathbb{R}}
 |\,\partial_xu (t, x)\,| 
\, \bigg)= 0.
$$
}

\medskip

The proofs of Theorem 1.1 
is given 
by a technical energy method. 

\medskip

This paper is organized as follows. 
In Section 2, We reformulate the problem 
in terms of the deviation from 
the asymptotic state. 
In order to show the asymptotics, 
we establish the {\it a priori} estimates 
by using the technical energy method in Section 3.  

\medskip

{\bf Some Notation.}\quad 
We denote by $C$ generic positive constants unless 
they need to be distinguished. 
In particular, use 
$C_{\alpha, \beta, \cdots }$ 
when we emphasize the dependency on $\alpha,\: \beta,\: \cdots $.

For function spaces, 
${L}^p = {L}^p(\mathbb{R})$ and ${H}^k = {H}^k(\mathbb{R})$ 
denote the usual Lebesgue space and 
$k$-th order Sobolev space on the whole space $\mathbb{R}$ 
with norms $||\cdot||_{{L}^p}$ and $||\cdot||_{{H}^k}$, 
respectively. 

\bigskip 

\noindent
\section{Reformulation of the problem} 
In this section, we reformulate our problem (1.1) 
in terms of the deviation from the asymptotic state. 
Now letting 
\begin{equation}
u(t, x) = \tilde{u} + \psi(t, x), 
\end{equation}
we reformulate the problem (1.1) in terms of 
the deviation $\psi $ from $\tilde{u}$ as 
\begin{eqnarray}
 \left\{\begin{array}{ll}
  \partial _t\psi 
  + \partial_x \big( \, f(\psi+\tilde{u}) \, \big) 
  =-\delta \, \partial_x^3 \psi - \nu \, \partial_x^4 \psi 
     \quad  \bigl( \, t>0,\: x\in \mathbb{R} \, \bigr), \\[5pt]
  \psi(0, x) = \psi_0(x) 
  := u_0(x)-\tilde{u} %
  \rightarrow 0 \quad (x\rightarrow \pm \infty).
 \end{array}
 \right.\,
\end{eqnarray}
Then we look for 
the unique global in time solution 
$\psi $ which has the asymptotic behavior 
\begin{equation}
\displaystyle{
\sup_{x \in \mathbb{R}}|\,  \psi(t, x) \, |
\to  0, \quad 
\sup_{x\in \mathbb{R}}
 |\,\partial_x\psi (t, x)\,|\to  0 
\qquad (t\to  \infty)}. 
\end{equation}
Here we note 
that $\psi_0 \in H^2$.
Then the corresponding theorems 
for $\psi$ to Theorems 1.1 we should prove is stated as follows. 

\medskip

\noindent
{\bf Theorem 2.1.} (Global Existence){\bf .}\quad{\it
Assume 
the convective flux $f \in C^2(\mathbb{R})$ 
satisfy {\rm(1.3)} and the initial data satisfy
$\psi_0 \in H^2$. 
Then the Cauchy problem {\rm(2.2)} has a 
unique global in time 
solution $\psi$ 
satisfying 
\begin{eqnarray*}
\left\{\begin{array}{ll}
\psi \in C^0\bigl( \, [\, 0, \, \infty ) \, ; H^2 \, \bigr), \\[5pt]
\partial _x^2 \psi \in L^2\bigl( \, 0,\, \infty \, ; H^2 \, \bigr),
\end{array} 
\right.\,
\end{eqnarray*}
and the asymptotic behavior 
$$
\lim _{t \to \infty}
\bigg( \,
\sup_{x\in \mathbb{R}}
|\,\psi (t, x)\,| 
+\sup_{x\in \mathbb{R}}
 |\,\partial_x\psi (t, x)\,| 
\, \bigg)= 0. 
$$
}

%
%

\medskip

To accomplish the proofs of Theorem 2.1, 
we prepare the local existence precisely, 
we formulate 
the problem (2.2) at general 
initial time $\tau \ge 0$: 
\begin{eqnarray}
 \left\{\begin{array}{ll}
  \partial _t\psi 
  + \partial_x \left( \, f(\psi+\tilde{u}) \, \right) 
  = - \delta \, \partial_x^3 \psi - \nu \, \partial_x^4 \psi 
     \quad  \bigl( \, t>\tau,\: x\in \mathbb{R} \, \bigr), \\[5pt]
  \psi(\tau, x) = \psi_\tau(x) 
  := u_\tau(x)-\tilde{u} %
  \rightarrow 0 \quad (x\rightarrow \pm \infty).
 \end{array}
 \right.\,
\end{eqnarray}

Then the local existence theorem is stated as follows 
(since the proof is standard, 
we state only here (cf. \cite{dua-fan-kim-xie}, \cite{dua-zha}, 
\cite{rua-gao-che}, \cite{wan-zhu})). 

\medskip

\noindent
{\bf Theorem 2.2} (Local Existence){\bf .}\quad{\it
For any $ M > 0 $, there exists a positive constant 
$t_0=t_0(M)$ not depending on $\tau$ 
such that if 
$\psi_{\tau} \in H^2$ and 
$$
\displaystyle{ 
\| \, \psi_{\tau} \, \|_{L^2} + \| \, \partial_x^2\psi_{\tau} \, \|_{L^2} \leq M}, 
$$
then the Cauchy problem {\rm (2.4)} has a unique solution $\psi$ 
on the time interval $[\, \tau, \, \tau+t_{0}(M)\, ]$ satisfying 
\begin{eqnarray*}
\left\{\begin{array}{ll}
\psi \in C^0
\bigl( \, [\, \tau, \, \tau+t_{0}\, ] \, ; H^2 \bigr),\\[5pt]
\partial _x^2 \psi  \in L^2\bigl( \, \tau, \, \tau+t_{0} \, ; H^2 \bigr),\\[5pt]
\displaystyle{
\sup_{t \in [\tau, \tau + t_{0}]} \, 
\left( \, 
\| \, \psi(t) \, \|_{L^2} + \| \, \partial_x^2\psi(t) \, \|_{L^2}
\, \right)
}
\leq 2 \, M. 
\end{array} 
\right.\,
\end{eqnarray*}
}

\medskip

Next, we state the {\it a priori} estimates as follows. 

\medskip

\noindent
{\bf Theorem 2.3} ({\it A Priori} Estimates){\bf .}\quad{\it 
Under the same assumptions as in Theorem 2.1, 
for any initial data 
$\psi_0 \in H^2$, 
there exists a positive constant $C_{\psi_0}$
such that 
if the Cauchy problem {\rm (2.1)}
has a solution 
$\phi$ 
on the time interval $[\, 0, \, T\, ]$ satisfying 
\begin{eqnarray*}
\left\{\begin{array}{ll}
\phi \in C^0
\bigl( \, [\, 0, \, T\, ] \, ; H^2 \bigr),\\[5pt]
\partial _x^2 \psi  \in L^2\bigl( \, 0, \, T \, ; H^2 \bigr),
\end{array} 
\right.\,
\end{eqnarray*}
for some positive constant $T$, 
then it holds that 
\begin{align}
\begin{aligned}
&\| \, \psi(t) \, \|_{H^2}^{2} 
  + \int^{t}_{0} 
\| \, \partial_{x}^2\psi(\tau) \, \|_{H^2}^{2} \, \mathrm{d}\tau
+\int^{t}_{0} 
\left( \, 
\sup_{x \in \mathbb{R}} |\, \psi(\tau,x) \, |
\, \right)^8
\, \mathrm{d}\tau \\
&
+ \int^{t}_{0} 
\left( \, 
\sup_{x \in \mathbb{R}} |\, \partial_x\psi(\tau,x) \, |
\, \right)^{\frac{8}{3}}
\, \mathrm{d}\tau 
\leq C_{\psi_0} 
\quad \big( \, t \in [\, 0, \, T\, ] \, \big).
\end{aligned}
\end{align}
}

\medskip

Combining the local existence Theorem 2.2 together with 
the each {\it a priori} estimates, Theorem 2.3 , we can obtain 
global existence Theorem 2.1. 
In fact, we can obtain the unique global in time solutions $\psi$ to (2.2) in 
Theorem 2.2 satisfying 
\begin{eqnarray*}
\left\{\begin{array}{ll}
\psi \in C^0\bigl( \, [\, 0, \, \infty ) \, ; H^2 \, \bigr), \\[5pt]
\partial _x^2 \psi \in L^2\bigl( \, 0,\, \infty \, ; H^2 \, \bigr),
\end{array} 
\right.\,
\end{eqnarray*}
and 
\begin{equation}
\sup_{t \ge 0}
  \| \, \psi(t) \, \|_{H^2}^{2} 
+ \int^{\infty}_{0} 
\| \, \partial_{x}^2\psi(t) \, \|_{H^2}^{2} \, \mathrm{d}t
< \infty 
\end{equation}
which yields 
\begin{equation}
\int _0^{\infty }\bigg|\,\frac{\mathrm{d}}{\mathrm{d}t} 
\|  \, \partial_x^2\psi(t)  \, \|_{L^2}^2 \,\bigg| \, \mathrm{d}t
< \infty.
\end{equation}
%
%
%
%
We immediately have from (2.6) and (2.7) that 
\begin{equation}
\|  \, \partial_x\psi(t)  \, \|_{L^2} \to 0\quad (t \to \infty).
\end{equation}
Further from (2.8) with $T\rightarrow \infty$ and 
\begin{equation}
\sup_{x \in \mathbb{R}}|\,\partial_x\psi(t,x)\,|
\leq \sqrt{2}\, \| \,\partial_x\psi(t)\, \|^{\frac{1}{2}}_{L^2} 
\| \,\partial_x^2\psi(t)\, \| ^{\frac{1}{2}}_{L^2} 
\quad \big( \, t \ge 0 \, \big), 
\end{equation} 
we obtain the asymptotic behavior (2.3). 

Thus Theorem 2.1 is proved.


\bigskip 

\noindent
\section{{\it A priori} estimates}
In this section, under the assumption 
\begin{equation}
|\, f''(u) \,| \leq O(1)\, \big( \, 1 + |\, u \,|^q \, \big) \quad \big( \, q\ge0 \, \big), 
\end{equation}
we show the following {\it a priori} estimate for $\psi$ in Theorem 2.3. 
To do that, we prepare the following basic estimate. 

\medskip

\noindent
{\bf Proposition 3.1.}\quad {\it
For $q \ge 0$, it follows that 
\begin{equation*}
\| \, \psi(t) \, \|_{L^2}^{2} +
2\, \nu \,
\int^{t}_{0} 
\| \, \partial_{x}^2\psi(\tau) \, \|_{L^2}^{2} \, \mathrm{d}\tau 
=\| \, \psi_0 \, \|_{L^2}^{2}
\quad \big( \, t \in [\, 0, \, T\, ] \, \big). 
\end{equation*}
}

\medskip

{\bf Proof of Proposition 3.1.}
Multiplying the equation in (2.2) by $\phi$ and integrating it with respect to $x$, 
we have, after integration by parts, 
\begin{equation}
\frac{1}{2} \, \frac{\mathrm{d}}{\mathrm{d}t} \, \| \, \psi(t) \, \|_{L^2}^2 
+ \nu \, \| \, \partial_{x}^2\psi(t) \, \|_{L^2}^{2} 
=0.
\end{equation}
Next, integrating (3.2) with respect to $t$, we immediately get the desired estimate.

Thus, we complete the proof of Proposition 3.1. 

\medskip

From Proposition 3.1, we have the next lemma. 

\medskip

\noindent
{\bf Lemma 3.2.}\quad {\it
There exists a positive constant $C_{\psi_0}$
such that 
\begin{equation*}
\int^{t}_{0} 
\left( \, 
\sup_{x \in \mathbb{R}} |\, \psi(\tau,x) \, |
\, \right)^8
\, \mathrm{d}\tau 
+ \int^{t}_{0} 
\left( \, 
\sup_{x \in \mathbb{R}} |\, \partial_x\psi(\tau,x) \, |
\, \right)^{\frac{8}{3}}
\, \mathrm{d}\tau 
\leq C_{\phi_0} 
\quad \big( \, t \in [\, 0, \, T\, ] \, \big). 
\end{equation*}
}

\medskip

{\bf Proof of Lemma 3.2.}
By using the Sobolev inequality and the integration by parts, 
we get 
\begin{align}
\begin{aligned}
\sup_{x \in \mathbb{R}}|\,\psi(t,x)\,|
&\le \sqrt{2}\, \| \,\psi(t)\, \|^{\frac{1}{2}}_{L^2} 
\| \,\partial_x\psi(t)\, \| ^{\frac{1}{2}}_{L^2}\\
&\le \sqrt{2}\, \| \,\psi(t)\, \|^{\frac{3}{4}}_{L^2} 
\| \,\partial_x^2\psi(t)\, \| ^{\frac{1}{4}}_{L^2}
\quad \big( \, t \in [\, 0, \, T\, ] \, \big),
\end{aligned}
\end{align}
\begin{align}
\begin{aligned}
\sup_{x \in \mathbb{R}}|\,\psi(t,x)\,|
&\le \sqrt{2}\, \| \,\partial_x\psi(t)\, \|^{\frac{1}{2}}_{L^2} 
\| \,\partial_x^2\psi(t)\, \| ^{\frac{1}{2}}_{L^2}\\
&\le \sqrt{2}\, \| \,\psi(t)\, \|^{\frac{1}{4}}_{L^2} 
\| \,\partial_x^2\psi(t)\, \| ^{\frac{3}{4}}_{L^2}
\quad \big( \, t \in [\, 0, \, T\, ] \, \big).
\end{aligned}
\end{align}
From (3.3) and (3.4), noting Proposition 3.1, we immediately have the desired estimate.

Thus, the proof is complete.

\medskip

Next, we state the {\it a priori} estimate for $\partial_{x}^2\psi$ as follows. 

\medskip

\noindent
{\bf Proposition 3.3.}\quad {\it
For $0\le q \le 5$, there exists a positive constant $C_{\psi_0}$
such that 
\begin{equation*}
\| \, \partial_{x}^2\psi(t) \, \|_{L^2}^{2} 
  + \int^{t}_{0} 
\| \, \partial_{x}^4\psi(\tau) \, \|_{L^2}^{2} \, \mathrm{d}\tau
\leq C_{\psi_0} 
\quad \big( \, t \in [\, 0, \, T\, ] \, \big).
\end{equation*}
}

Once Proposition 3.3 holds true, by using Proposition 3.1, 
we can estimate as follows. 
\begin{equation}
\| \, \partial_{x}\psi \, \|_{L^2}^{2} 
\le \| \, \psi \, \|_{L^2} \, \| \, \partial_{x}^2\psi \, \|_{L^2} 
\le C_{\phi_0},
\end{equation}
\begin{equation}
\int^{t}_{0} 
\| \, \partial_{x}^3\psi \, \|_{L^2}^{2} \, \mathrm{d}\tau
\le \left( \, \int^{t}_{0} 
    \| \, \partial_{x}^2\psi \, \|_{L^2}^{2} \, \mathrm{d}\tau \, \right)^{\frac{1}{2}} \, 
    \left( \, \int^{t}_{0} 
    \| \, \partial_{x}^4\psi \, \|_{L^2}^{2} \, \mathrm{d}\tau \, \right)^{\frac{1}{2}}
\le C_{\phi_0}, 
\end{equation}
for $t \in [\,0, \, T \,]$. 
From the uniform estimates (3.5) and (3.6), 
we immediately get the {\it a priori} estimate for $\partial_{x}\psi$ as follows. 

\medskip

\noindent
{\bf Proposition 3.4.}\quad {\it
For $0\le q \le 5$, there exists a positive constant $C_{\psi_0}$
such that 
\begin{equation*}
\| \, \partial_{x}\psi(t) \, \|_{L^2}^{2} 
  + \int^{t}_{0} 
\| \, \partial_{x}^3\psi(\tau) \, \|_{L^2}^{2} \, \mathrm{d}\tau
\, \mathrm{d}x\mathrm{d}\tau 
\leq C_{\psi_0} 
\quad \big( \, t \in [\, 0, \, T\, ] \, \big).
\end{equation*}
}

From Propositions 3.1, 3.3-3.4, by using the Sobolev inequality, we have 
the following uniform boundedness of $\psi$ and $\partial_{x}\psi$ as follows. 

\medskip

\noindent
{\bf Lemma 3.5.}\quad {\it
There exists a positive constant $C_{\psi_0}$ such that 
$$
\sup_{x \in \mathbb{R}} \, | \, \psi(t,x) \, |\le C_{\psi_0}, 
\quad \sup_{x \in \mathbb{R}} \, | \, \partial_{x}\psi(t,x) \, |\le C_{\psi_0} 
\quad \big( \, t \in [\, 0, \, T\, ] \, \big).
$$
}


\medskip
By combining Propositions 3.1, 3.3-3.4 and Lemma 3.2, 
we can obtain Theorem 2.3. 
Therefore, in order to complete the proof of Theorem 2.3, we finally prove Proposition 3.3. 

\medskip

{\bf Proof of Proposition 3.3.}
Multiplying the equation in (2.2) by 
$\partial_x^4 \psi$, and integrating the resultant formula 
with respect to $x$, we have 
\begin{equation}
\frac{1}{2} \, \frac{\mathrm{d}}{\mathrm{d}t} \, \| \, \partial_x\psi(t) \, \|_{L^2}^2 
+  \nu \, \| \, \partial_x^4\psi(t) \, \|_{L^2}^{2}  
= - \int^{\infty}_{-\infty} 
   \partial_x^4\psi \, 
   \partial_x
   \bigl( \, f(\psi +\tilde{u}) \, \bigr) 
   \, \mathrm{d}x.
\end{equation}
The right-hand side of (3.7) becomes 
\begin{align}
\begin{aligned}
&- \int^{\infty}_{-\infty} 
   \partial_x^4\psi \, 
   \partial_x
   \bigl( \, f(\psi +\tilde{u}) \, \bigr) 
   \, \mathrm{d}x \\
&=\int^{\infty}_{-\infty} 
   f''(\psi + \tilde{u}) \, | \, \partial_x\psi \, |^2 \, \partial_x^3 \psi
   \, \mathrm{d}x+ \int^{\infty}_{-\infty} 
   f'(\psi +\tilde{u}) \, \partial_x^2 \psi \, \partial_x^3 \psi
   \, \mathrm{d}x\\
&=\int^{\infty}_{-\infty} 
   f''(\psi + \tilde{u}) \, | \, \partial_x\psi \, |^2 \, \partial_x^3 \psi
   \, \mathrm{d}x
   -\frac{1}{2} \, \int^{\infty}_{-\infty} 
   f''(\psi +\tilde{u}) \, \partial_x \psi \, | \, \partial_x^2 \psi \,|^2
   \, \mathrm{d}x.
\end{aligned}
\end{align}
From (3.1), 
by making use of the Cauchy-Schwarz, Sobolev and Young inequalities, and 
integration by parts, 
we can estimate the first term on the right-hand side of (3.8) as follows. 
\begin{align}
\begin{aligned}
&
\left| \, 
\int^{\infty}_{-\infty} 
f''(\psi + \tilde{u}) \, | \, \partial_x\psi \, |^2 \, \partial_x^3 \psi
   \, \mathrm{d}x 
   \, \right|\\
&\le C_q \, \int^{\infty}_{-\infty} 
| \, \partial_x\psi \, |^2 \, | \, \partial_x^3 \psi \, |
   \, \mathrm{d}x
+ C_q \, \int^{\infty}_{-\infty} 
     | \, \phi \, |^q \, | \, \partial_x\psi \, |^2 \, | \, \partial_x^3 \psi \, |
   \, \mathrm{d}x, 
\end{aligned}
\end{align}
\begin{align}
\begin{aligned}
&C_q \, \int^{\infty}_{-\infty} 
     | \, \partial_x\psi \, |^2 \, | \, \partial_x^3 \psi \, |
   \, \mathrm{d}x\\
&\le C_q \, \| \, \partial_x\psi \, \|_{L^2}^2 
     \, \| \, \partial_x^3\psi \, \|_{L^2}^{\frac{1}{2}}
     \, \| \, \partial_x^4\psi \, \|_{L^2}^{\frac{1}{2}}\\
&\le C_q \, \| \, \psi \, \|_{L^2} 
     \, \| \, \partial_x^2\psi \, \|_{L^2}^{\frac{5}{4}}
     \, \| \, \partial_x^4\psi \, \|_{L^2}^{\frac{3}{4}}\\
&\le \frac{\nu}{8} \, \| \, \partial_x^4\psi \, \|_{L^2}^2 
     + C_{q, \nu} \, \| \, \psi \, \|_{L^2}^{\frac{8}{5}} 
       \, \| \, \partial_x^2\psi \, \|_{L^2}^2, 
\end{aligned}
\end{align}
\begin{align}
\begin{aligned}
&C_q \, \int^{\infty}_{-\infty} 
     | \, \psi \, |^q \, | \, \partial_x\psi \, |^2 \, | \, \partial_x^3 \psi \, |
   \, \mathrm{d}x\\
&\le C_q \, \| \, \psi \, \|_{L^{\infty}}^{q}
     \, \| \, \partial_x\psi \, \|_{L^2}^2
     \, \| \, \partial_x^3\psi \, \|_{L^{\infty}}\\
&\le C_q \, \| \, \psi \, \|_{L^{\infty}}^{q}
     \, \| \, \psi \, \|_{L^2}
     \, \| \, \partial_x^2\psi \, \|_{L^2}^{\frac{5}{4}}
     \, \| \, \partial_x^4\psi \, \|_{L^2}^{\frac{3}{4}}\\
&\le \frac{\nu}{8} \, \| \, \partial_x^4\psi \, \|_{L^2}^2 
+ C_{q, \nu} \, \| \, \psi \, \|_{L^{\infty}}^{\frac{8q}{5}}
     \, \| \, \psi \, \|_{L^2}^{\frac{8}{5}}
     \, \| \, \partial_x^2\psi \, \|_{L^2}^2\\
&\le \frac{\nu}{8} \, \| \, \partial_x^4\psi \, \|_{L^2}^2 
+ C_{q, \nu} \, \, \left( \, 1+ \| \, \psi \, \|_{L^{\infty}}^8 \, \right)
     \| \, \psi \, \|_{L^2}^{\frac{8}{5}}
     \, \| \, \partial_x^2\psi \, \|_{L^2}^2
     \quad \big( \, 0 \le q \le 5 \, \big).
\end{aligned}
\end{align}
Similarly, the second term on the right-hand side of (3.8) can be estimated as follows. 
\begin{align}
\begin{aligned}
&\left| \, 
\frac{1}{2} \, \int^{\infty}_{-\infty} 
   f''(\psi +\tilde{u}) \, \partial_x \psi \, | \, \partial_x^2 \psi \,|^2
   \, \mathrm{d}x
\, \right|\\
&\le C_q \, \int^{\infty}_{-\infty} 
     | \, \partial_x \psi \, | \, | \, \partial_x^2 \psi \, |^2
   \, \mathrm{d}x
+ C_q \, \int^{\infty}_{-\infty} 
     | \, \psi \, |^{q} \, | \, \partial_x \psi \, | \, | \, \partial_x^2 \psi \, |^2
   \, \mathrm{d}x,
\end{aligned}
\end{align}
\begin{align}
\begin{aligned}
C_q \, \int^{\infty}_{-\infty} 
     | \, \partial_x \psi \, | \, | \, \partial_x^2 \psi \, |^2
   \, \mathrm{d}x
&\le C_q \, \| \, \partial_x\psi \, \|_{L^2}^{\frac{1}{2}} 
     \, \| \, \partial_x^2\psi \, \|_{L^2}^{\frac{5}{2}}\\
&\le C_q \, \| \, \psi \, \|_{L^2}^{\frac{1}{4}} 
     \, \| \, \partial_x^2\psi \, \|_{L^2}^{\frac{11}{4}}\\
&\le C_q \, \| \, \psi \, \|_{L^2}^{\frac{1}{4}} 
     \, \| \, \partial_x^2\psi \, \|_{L^2}^2
     \, \left( \, 1+ \| \, \partial_x^2\psi \, \|_{L^2}^2 \, \right), 
\end{aligned}
\end{align}
\begin{align}
\begin{aligned}
&
C_q \, \int^{\infty}_{-\infty} 
     | \, \psi \, |^{q} \, | \, \partial_x \psi \, | \, | \, \partial_x^2 \psi \, |^2
   \, \mathrm{d}x\\
&\le C_q \, \| \, \psi \, \|_{L^{\infty}}^{q} 
     \, \| \, \partial_x\psi \, \|_{L^{\infty}} 
     \, \| \, \partial_x^2\psi \, \|_{L^2}^{2}\\
&\le C_q \, \left( \, \| \, \psi \, \|_{L^{\infty}}^{\frac{8q}{5}}
     + \| \, \partial_x\psi \, \|_{L^{\infty}}^{\frac{8}{3}} \, \right)
     \, \| \, \partial_x^2\psi \, \|_{L^2}^{2}\\
&\le C_q \, \left( \, 1+ \| \, \psi \, \|_{L^{\infty}}^8
     + \| \, \partial_x\psi \, \|_{L^{\infty}}^{\frac{8}{3}} \, \right)
     \, \| \, \partial_x^2\psi \, \|_{L^2}^{2}
     \quad \big( \, 0 \le q \le 5 \, \big).
\end{aligned}
\end{align}
Noting Lemma 3.2, 
substituting (3.8)-(3.14) into (3.7), 
integrating the resultant formula with respect to $t$ 
and further using the Gronwall inequality, 
we obtain the desired formula, Proposition 3.3. 

Thus, we complete the proof of Theorem 2.3 from Propositions 3.1, 3.3-3.4. 

\bigskip









\bibliographystyle{model6-num-names}
\bibliography{<your-bib-database>}







\end{document}